\documentclass[11pt]{article}
\usepackage{amssymb,amsmath,amsthm}
\newtheorem{theorem}{Theorem}

\newtheorem{lemma}{Lemma}

\newtheorem{remark}{Remark}
\date{}
\numberwithin{equation}{section} \numberwithin{theorem}{section}
\numberwithin{lemma}{section} \numberwithin{corollary}{section}
\numberwithin{remark}{section} \numberwithin{proposition}{section}
\numberwithin{definition}{section}
\usepackage{color}
\usepackage{curves}
\usepackage[dvips]{graphics}
\usepackage{amsfonts}
\usepackage{amsmath}
\usepackage{amssymb}
\usepackage[latin1]{inputenc}
\usepackage[normalem]{ulem}
\usepackage{t1enc}
\usepackage{fancybox}
\usepackage{epsfig}
\usepackage{minitoc}
\usepackage{fancyhdr}

\newcommand{\n}{\noindent}

\newcommand{\vs}{\vskip}
\begin{document}

\title{Uniqueness of Solution of the Unsteady Filtration Problem
in Heterogeneous Porous Media }

\author{A. Lyaghfouri$^1$ and E. Zaouche$^2$\\
$^1$ American University of Ras Al Khaimah\\
Ras Al Khaimah, UAE\\
$^2$ Ecole Normale Sup\'{e}rieure \\
Algiers, Algeria}
\maketitle

\begin{abstract}
We establish uniqueness of the solution of the unsteady state dam problem
in the heterogeneous and rectangular case assuming the dam wet at the
bottom and dry near to the top.
\end{abstract}

\begin{flushleft}
\end{flushleft}

\begin{flushleft}
\end{flushleft}

\section{Introduction}\label{s1}
We consider a heterogeneous porous medium supplied by several reservoirs
of a fluid, represented by a bounded domain $\Omega$ of ${\mathbb R}^n$ with locally Lipschitz boundary
 $\partial\Omega=\Gamma_1\cup \Gamma_2$, where $\Gamma_1$ is the impervious
 part of the boundary,
 $\Gamma_2$ is the part in contact with either air or the fluid reservoirs.

\n The fluid infiltrates through $\Omega$
obeying to Darcy's law
$$v = -a(x).\nabla (u + x_n),$$
where $a(x)=(a_{ij}(x))_{ij}$ is the $n\times n$ permeability matrix
of the medium,
$x=(x_1,...,x_n)$, $v$ is the fluid velocity and $u$ its pressure.

\n We are concerned with the problem of finding the pressure $u$ and the
 saturation $\chi$ of the fluid inside $\Omega$. Using the mass conservation law,  
Darcy's law, the boundary conditions and the initial data, we obtain the following strong
formulation for our problem (see~\cite{[CL]}):
\begin{equation*}{\bf (SF)} \quad \left\{
\begin{aligned}
\quad u \geq 0, \; 0\leq \chi\leq 1,\; u(1-\chi ) &= 0
&\quad & \text{ in }Q \\
div (a(x)(\nabla u+\chi e))-(\alpha u+\chi)_t &= 0 & \quad & \text{ in } Q \\
u &=\phi &\quad & \text{ on } \Sigma_2 \\
(\alpha u+\chi)(\cdot,0)&=\alpha u_0+\chi_0 & \quad & \text{ in } \Omega \\
a(x)(\nabla u+\chi e)\cdot \nu &=0  & \quad
& \text{ on } \Sigma_1 \\
a(x)(\nabla u+\chi e)\cdot \nu &\leq 0  & \quad
& \text{ on } \Sigma_4
\end{aligned} \right. \hspace{3.3cm}
\end{equation*}

\n  where $\alpha, T$ are  positive numbers, $Q=\Omega\times (0,T)$,
$\Sigma_1=\Gamma_1\times(0,T)$ is the impervious part of $\partial\Omega$, 
$\Sigma_2=\Gamma_2\times(0,T)$ is the pervious part,
$\Sigma_3=\Sigma_2\cap \{\phi>0\}$ is the part covered by fluid,
and $\Sigma_4=\Sigma_4\cap \{\phi=0\}$ is the part where the
fluid flows outside $\Omega$.

\n $\phi$ is a nonnegative Lipschitz continuous function defined in  $\overline{Q},$
$\nu$ is the outward unit normal vector to $\partial\Omega$, $e=(0,...,0,1)\in {\mathbb R}^n$,
$u_0, \chi_0: \Omega \longrightarrow \mathbb{R}$ are functions satisfying for a positive constant
$M$
\begin{eqnarray}\label{e1.1}
0\leq u_0(x)\leq M,\quad 0\leq \chi_0(x)\leq 1\quad\text{for a.e. } x\in\Omega.
\end{eqnarray}

\n For $a(x)$, we assume that we have for two positive constants $\lambda$ and $\Lambda$
\begin{eqnarray}\label{e1.3-1.4}
&&\forall \xi\in \mathbb{R}^n,\quad\text{for a.e. } x\in\Omega  \quad \lambda |\xi|^2\leq
a(x).\xi.\xi.\\
&&\forall \xi\in \mathbb{R}^n,\quad\text{for a.e. } x\in\Omega  \quad  |a(x).\xi|\leq \Lambda|\xi|.
\end{eqnarray}
Moreover, we assume that
\begin{eqnarray}\label{e1.5}
&&div(a(x)e) \in L^{2}(\Omega).
\end{eqnarray}

\n Using the strong formulation, we are led to the following weak formulation
\begin{equation*}{\bf (P)} \quad \left\{
\begin{aligned}
&\text{Find }
(u, \chi) \in L^2(0,T;H^1(\Omega))\times L^\infty (Q) \text{ such that}:\\[0.2cm]
(i)&\quad u \geq 0, \; 0\leq \chi\leq 1,\; u.(1-\chi) = 0  \quad \text{ a.e. in }Q \\[0.1cm]
(ii)&\quad u=\phi \hspace{0.34cm} \text{ on } \Sigma_2 \\
(iii)&\quad \displaystyle{\int_Q \big[ a(x)\big(\nabla u +\chi e\big)
  \cdot\nabla\xi-(\alpha u+\chi)\xi_t\big] dx\,dt}\\
&\qquad\qquad  \leq \displaystyle{\int_\Omega (\chi_0+\alpha u_0)\xi(x,0)\,dx}  \\
&\quad \forall \xi \in H^1(Q),\; \xi=0 \text{ on }
   \Sigma_3,\; \xi \geq 0 \text{ on } \Sigma_2,\;\\
& \quad\xi(x,T)=0 \,\text{ for a.e. } x\in \Omega.
\end{aligned} \right.
\end{equation*}

\vs0.3cm\n For the existence of a solution of the problem $(P)$ in the homogeneous case
$(a(x)=I_n)$, we refer to \cite{[CG]} and  \cite{[C1]} respectively in the incompressible $(\alpha=0)$
and compressible $(\alpha>0)$ cases. For the heterogeneous case, we refer to \cite{[Z]}
in a more general framework under assumptions (1.1)-(1.4) for both incompressible
and compressible cases. For the incompressible case with nonlinear Darcy's law,
we refer to \cite{[L2]}, \cite{[L3]} and \cite{[CL]} respectively for
Dirichlet, Neuman and generalized boundary conditions.
Regarding regularity of the solution, we refer to \cite{[CG]} and \cite{[C1]},
where it has been proved when $a(x)=I_n$ that
$\chi \in C^0([0,T];L^p(\Omega))$ for all $p\geq 1$ in both incompressible and
compressible cases, and that $u \in C^0([0,T];L^p(\Omega))$ for all $ 1\leq p\leq 2$
in the compressible case. Extensions to the quasilinear and incompressible case
were obtained in \cite{[L2]}, \cite{[L3]} and \cite{[L4]} in both homogeneous
and nonhomogeneous frameworks. The authors of this paper were recently able to extend the
above regularity result in \cite{[LZ]} to a more general framework under weaker assumptions
on the data.

\vs 0.2cm\n In this paper, we are mainly concerned with the uniqueness of the solution
of the problem $(P)$. This question was first addressed for a rectangular homogeneous dam
in \cite{[T]} and \cite{[DF]} respectively for a formulation based on quasi-variational
inequalities and for the formulation $(P)$ with a dam wet at the bottom and dry near
to the top in the second case.
Uniqueness of the solution for a homogeneous dam with general geometry was established
by the method of doubling variables in \cite{[C1]}, but it is not obvious wether it works
in the heterogeneous situation.
Extensions to a quasilinear operator modeling incompressible fluid flow governed
by a nonlinear Darcy's law with Dirichlet, or Neuman boundary conditions were obtained in
\cite{[L2]} and \cite{[L3]} respectively.

\vs 0.2cm\n Our main result in this work is the uniqueness of problem $(P)$ solution
for a heterogeneous and rectangular porous medium assuming it remains 
wet near the bottom and dry close to the top. Our method is inspired by an idea from \cite{[DF]}
in the homogeneous case and relies on solution regularity that has been recently obtained in \cite{[LZ]}.
Our uniqueness result is new in the heterogeneous and rectangular framework, but most likely the technic 
is limited to that particular geometry like shape.

\section{Preliminary Results}\label{s2}

\n In this work we shall be interested in the following situation of a two
dimensional rectangular dam $\Omega=(0,L)\times (0,K)$, with $L,K>0$, 
and $\Gamma_{1}=[0,L]\times\{0\}$, $\Gamma_{2}=(\{0\}\times [0,K])\cup ([0,L]\times\{K\})\cup (\{L\}\times [0,K])$
(see Figure 1).

\n We also assume that
\begin{eqnarray}\label{e2.1-3}
&&a(x)e\in C^{0,1}(\overline{\Omega}),\\
&&div(a(x)e) \geq 0 \quad\text{a.e. in } \Omega,\\
&&\phi_0\leq\phi\leq \phi_1  \;\; \text{ on }\Sigma_{2},
\end{eqnarray}
\n where $\phi_0$ and $\phi_1$ are two nonnegative Lipschitz
continuous functions defined on $\overline{\Omega}$ and satisfying for some
$\epsilon_{0}>0$ small enough
\begin{figure}[t]
\setlength{\unitlength}{0.5mm}
\begin{center}

\begin{picture}(130,70)(-20,-20)
\put(2,-22){\line(1,0){98}}
\put(26,-22){\line(0,1){60}}
\put(76,-22){\line(0,1){60}}
\put(26,38){\line(1,0){50}}

\put(45,40){$\Gamma_2$}

\put(0,23){$y=h_1$}

\put(2,19){\line(1,0){24}}
\textcolor{blue}{\multiput(2,18.5)(1,0){24}{.}
\multiput(2,17.5)(1,0){24}{.}
\multiput(2,16.5)(1,0){24}{.}
\multiput(2,15.5)(1,0){24}{.}
\multiput(2,14.5)(1,0){24}{.}
\multiput(2,13.5)(1,0){24}{.}
\multiput(2,12.5)(1,0){24}{.}
\multiput(2,11.5)(1,0){24}{.}
\multiput(2,10.5)(1,0){24}{.}
\multiput(2,9.5)(1,0){24}{.}\multiput(2,8.5)(1,0){24}{.}
\multiput(2,7.5)(1,0){24}{.}\multiput(2,6.5)(1,0){24}{.}
\multiput(2,5.5)(1,0){24}{.}\multiput(2,4.5)(1,0){24}{.}
\multiput(2,3.5)(1,0){24}{.}\multiput(2,2.5)(1,0){24}{.}
\multiput(2,1.5)(1,0){24}{.} \multiput(2,0.5)(1,0){24}{.}
\multiput(2,-0.5)(1,0){24}{.}\multiput(2,-1.5)(1,0){24}{.}
\multiput(2,-2.5)(1,0){24}{.}\multiput(2,-3.5)(1,0){24}{.}
\multiput(2,-4.5)(1,0){24}{.}\multiput(2,-5.5)(1,0){24}{.}
\multiput(2,-6.5)(1,0){24}{.}\multiput(2,-7.5)(1,0){24}{.}
\multiput(2,-8.5)(1,0){24}{.}\multiput(2,-9.5)(1,0){24}{.}
\multiput(2,-10.5)(1,0){24}{.}\multiput(2,-11.5)(1,0){24}{.}
\multiput(2,-12.5)(1,0){24}{.}\multiput(2,-13.5)(1,0){24}{.}
\multiput(2,-14.5)(1,0){24}{.}\multiput(2,-15.5)(1,0){24}{.}
\multiput(2,-16.5)(1,0){24}{.}\multiput(2,-17.5)(1,0){24}{.}
\multiput(2,-18.5)(1,0){24}{.}\multiput(2,-19.5)(1,0){24}{.}
\multiput(2,-20.5)(1,0){24}{.}\multiput(2,-21.5)(1,0){24}{.}
\multiput(2,-22)(1,0){24}{.} \put(2,-22){\line(1,0){24}}}

\put(85,7.5){$y=h_2$}
\textcolor{blue}{\put(76,5){\line(1,0){24}}
\multiput(76,4)(1,0){24}{.}\multiput(76,3)(1,0){24}{.}
\multiput(76,2)(1,0){24}{.}\multiput(76,1)(1,0){24}{.}
\multiput(76,0)(1,0){24}{.}\multiput(76,-1)(1,0){24}{.}
\multiput(76,-2)(1,0){24}{.}\multiput(76,-3)(1,0){24}{.}
\multiput(76,-4)(1,0){24}{.}\multiput(76,-5)(1,0){24}{.}
\multiput(76,-6)(1,0){24}{.}
\multiput(76,-7)(1,0){24}{.}
\multiput(76,-8.5)(1,0){24}{.}
\multiput(76,-9.5)(1,0){24}{.}\multiput(76,-10.5)(1,0){24}{.}
\multiput(76,-11.5)(1,0){24}{.}\multiput(76,-12.5)(1,0){24}{.}
\multiput(76,-13.5)(1,0){24}{.}\multiput(76,-14.5)(1,0){24}{.}
\multiput(76,-15.5)(1,0){24}{.}\multiput(76,-15.5)(1,0){24}{.}
\multiput(76,-16.5)(1,0){24}{.}\multiput(76,-17.5)(1,0){24}{.}
\multiput(76,-18.5)(1,0){24}{.}\multiput(76,-19.5)(1,0){24}{.}
\multiput(76,-20.5)(1,0){24}{.}\multiput(76,-21.5)(1,0){24}{.}
}
\multiput(2,-24)(1,0){98}{-}

\put(45,-29){$\Gamma_1$}
 \put(28,-10){$\Gamma_2$}
  \put(65,-10){$\Gamma_2$}
\textcolor{blue}{\qbezier(26,19)(50,16)(76,8)}

\put(40,20){$\Gamma$}

\put(45,-40){Figure 1}
\end{picture}
\end{center}

\end{figure}
\vs 0.3cm
\begin{equation}\label{e2.4}
\left\{
\begin{aligned}
&\phi_0(0,x_2)=\phi_0(L,x_2)=(\epsilon_{0}-x_2)^{+} \\
&\phi_1(0,x_2)=\phi_1(L,x_2)=(K-\epsilon_{0}-x_2)^+\\
&~\phi_0(x_1,K)=\phi_1(x_1,K)=0.
\end{aligned} \right.
\end{equation}

\n Let us now denote by $(v_i,\gamma_i)$  the solution of the stationary problem corresponding
to $\phi_i,~~ i=0,1$ (see \cite{[ChL2]})

\begin{equation*}{\bf (P^{s}_{i})} \; \left\{
\begin{aligned}
\text{Find }& (v_i,\gamma_i) \in H^1(\Omega)\times L^\infty(\Omega)\text{ such that :  } \\
(i)&\quad v_i \geq 0, \; 0\leq \gamma_i\leq 1,\; v_i.(1-\gamma_i) = 0  \quad \text{ a.e. in }\Omega \\[0.1cm]
(ii)&\quad v_i=\phi_i \hspace{0.34cm} \text{ on } \Gamma_2 \\
(iii) &\int_{\Omega} a(x)(\nabla v_{i}+\gamma_i e).\nabla\xi dx\leq 0\\
& \forall \xi\in H^{1}(\Omega), \quad \xi=0~ \text{ on } ~\Gamma_{2}\cap\{\phi_{i}>0\},~~ \xi\geq 0 ~\text{ on }~  \Gamma_{2}\cap\{\phi_{i}=0\}.
\end{aligned} \right.
\end{equation*}

\vs 0.3cm\n We recall the following uniqueness result.

\begin{theorem}\label{t2.1} Assume that (2.1)-(2.2) hold.
Then the solution $(v_i,\gamma_i)$ of $(P^s_i)$ is unique and satisfies
\begin{equation}\label{2.5}
\gamma_i=\chi_{\{v_{i}>0\}}.
\end{equation}
\end{theorem}

\n \emph{Proof.} When $a(x)e\in C^1(\overline{\Omega})$, we refer to \cite{[ChL1]},
since we cannot have pools in a rectangular dam. One can also argue as in \cite{[ChL3]}.
When $a(x)e\in C^{0,1}(\overline{\Omega})$, one may combine Theorem 5.1 of \cite{[ChL2]}
and the proof of Theorem 6.3 in \cite{[ChL1]} to establish the result.

\qed

\vs 0,3cm \n The following properties for the solution of $(P^s_i)$ $(i=0,1)$ hold.

\begin{theorem}\label{t2.2} Assume that (2.1)-(2.4) hold.
Then we have
\begin{eqnarray}\label{e2.6-2.8}
&&(\epsilon_0-x_2)^+\leq v_0\leq v_1\leq (K-\epsilon_{0}-x_2)^+~~\text{ a.e. in }~ \Omega\\
&&\gamma_0=1~~ \text{ a.e. in }~ \Omega\cap\{0<x_2<\epsilon_0\}\\
&&\gamma_1=0 ~~ \text{ a.e. in }~ \Omega\cap\{K-\epsilon_0<x_2<K\}.
\end{eqnarray}
\end{theorem}

\vs 0.3cm\n \emph{Proof.} First, we remark that for each $k\in(0,K)$, $(k-x_2)^+$
satisfies the equation
$$\int_{\Omega} a(x)(\nabla (k-x_2)^++\chi_{\{(k-x_2)^+>0\}} e).\nabla\xi dx=\int_{\Omega\cap\{x_2<k\}} a(x)(-e+ e).\nabla\xi dx=0.$$
It follows that $((k-x_2)^+,\chi_{\{(k-x_2)^+>0\}})$ is a solution of the stationary dam problem for
the boundary Dirichlet data $(k-x_2)^+$ on $\Gamma_2$.
Adapting the proof of uniqueness in \cite{[ChL1]} and arguing as in \cite{[L1]},
and using (2.1)-(2.4), we obtain
\begin{eqnarray*}
&&(\epsilon_0-x_2)^+\leq v_0\leq v_1\leq (K-\epsilon_{0}-x_2)^+~~\text{ in }~ \Omega\\
&&\chi_{\{(\epsilon_0-x_2)^+>0\}}\leq \gamma_0\leq \gamma_1\leq \chi_{\{(K-\epsilon_0-x_2)^+>0\}}~~\text{ a.e. in }~ \Omega.
\end{eqnarray*}
Hence (2.6)-(2.8) follow. \qed

\vs 0.3cm\n

\begin{remark}\label{r2.1}
Theorems 2.1 and 2.2 remain true without the regularity assumption (2.1)
provided the following assumptions on the permeability matrix hold (see \cite{[L1]})
\begin{eqnarray*}
&& a_{12}=0\quad\text{a.e. in } \Omega,\\
&&{{\partial a_{22}}\over{\partial x_2}} \geq 0 \quad\text{in } \mathcal{D}'(\Omega).
\end{eqnarray*}
\end{remark}

\vs 0.3cm\n Next we will construct a solution corresponding
to a dam that is wet up to $x_2=\epsilon_0$ and dry above $x_2=K-\epsilon_0$ over
the whole interval $[0,T]$.

\begin{lemma}\label{l2.1} Assume that (2.1)-(2.4) hold and the initial data satisfies
\begin{eqnarray}\label{e2.9-2.10}
&&v_0\leq u_0\leq v_1\;\; \text{ a.e. in } \Omega.\\
&&\gamma_0\leq\chi_0\leq \gamma_1\;\; \text{ a.e. in } \Omega.
\end{eqnarray}
Then there exists a solution $(u,\chi)$ of problem (P) such that
\begin{eqnarray}\label{e2.11-2.12}
&&v_0\leq u\leq v_1\;\; \text{ a.e. in } Q\\
&&\gamma_0\leq \chi\leq \gamma_1\;\; \text{ a.e. in } Q.
\end{eqnarray}
\end{lemma}
\n \emph{Proof.} Let $v_{i\epsilon}$  be the solution
of the approximating problem of the stationary problem
$(P^{s}_{i}),~i=0,1$
\begin{equation*}{\bf (P^{s}_{i\epsilon})} \; \left\{
\begin{aligned}
 &\text{ Find } v_{i\epsilon} \in H^{1}(\Omega) \text{ such that :  } \\
  & (i)\quad v_{i\epsilon}=\phi_{i}  \quad \text{ on } \Gamma_{2}\\
 &(ii)\quad\int_{\Omega} a(x)(\nabla v_{i\epsilon}+H_{\epsilon}(v_{i\epsilon})e).\nabla\xi dx=0\\
&\quad \forall \xi\in H^{1}(\Omega),~~ \xi=0 \quad \text{ on }
\Gamma_{2},
\end{aligned} \right.
\end{equation*}
where $H_{\epsilon}(s)=\min(1,s^+/\epsilon)$ is an approximation of the Heaviside graph
$H(s)=[0,1]\chi_{\{0\}}+\chi_{(0,\infty)}.$

\n Let $u_{\epsilon}$ be the solution of the following approximating problem of the problem $(P)$
\begin{equation*} {\bf (P_{\epsilon})} \;\left\{
\begin{aligned}
&\text{ Find } u_{\epsilon}\in H^{1}(Q) \text{ such that}: \;\; \; \\
&(i)\quad u_{\epsilon}=\phi  \text{ on } \Sigma_2 \\
& (ii)\quad \displaystyle{\int_Q \big[ a(x)\big(\nabla
u_{\epsilon} +H_{\epsilon}(u_{\epsilon})e\big)
  \cdot\nabla\xi+\epsilon u_{\epsilon t}\xi_{t}-G_{\epsilon}(u_{\epsilon})\xi_t\big] dxdt}\\
  &\quad \quad+\int_{\Omega}G_{\epsilon}(u_{\epsilon}(x,T))\xi(x,T)dx=\int_{\Omega}(\alpha u_{0\epsilon}(x)+\chi_{0\epsilon}(x))\xi(x,0)dxdy\\
 &\quad \forall \xi \in H^1(Q),\; \xi=0 \text{ on }
   \Sigma_2,
\end{aligned} \right.
\end{equation*}
where $u_{0\epsilon}=\min(u_{0},v_{1\epsilon})$ and
$\chi_{0\epsilon}=\min(\chi_{0},H_{\epsilon}(v_{1\epsilon})).$\\
If $\xi \in H^1(Q),\; \xi=0 \text{ on }
   \Sigma_2,$ we have from $(P^{s}_{1\epsilon}) (ii):$
\begin{eqnarray}\label{e2.13}
&&\displaystyle{\int_Q \big[ a(x)\big(\nabla v_{1\epsilon}+H_{\epsilon}(v_{1\epsilon})e\big)
  \cdot\nabla\xi-(\alpha v_{1\epsilon}+H_{\epsilon}(v_{1\epsilon}))\xi_t\big] dxdt}\nonumber\\
  && +\int_{\Omega}(\alpha v_{1\epsilon}+H_{\epsilon}(v_{1\epsilon}))\xi(x,T)dx=\int_{\Omega}(\alpha
  v_{1\epsilon}+H_{\epsilon}(v_{1\epsilon}))\xi(x,0)dx.\nonumber\\
\end{eqnarray}
For $\delta>0,$ the function $\displaystyle{\xi_{\delta}=\frac{(u_{\epsilon}-v_{1\epsilon}-\delta)^{+}}{u_{\epsilon}-v_{1\epsilon}}}$
belongs to $ H^{1}(Q)$ and satisfies $ \xi_{\delta}=0$ on $\Sigma_{2}$ since $ \phi\leq\phi_{1}$ on
$\Sigma_{2}.$ Writing
(2.13) and $(P_{\epsilon})(ii)$ for $\xi=\xi_{\delta}$ and subtracting
the two identities from each other, we get by taking into account (2.9)-(2.10)
\begin{eqnarray}\label{e2.14}
&&\int_{Q}\Big[a(x)\big(\nabla (u_{\epsilon}-v_{1\epsilon})
+(H_{\epsilon}(u_{\epsilon})-H_{\epsilon}(v_{1\epsilon}))e\big).\nabla\xi_{\delta}+\epsilon(u_{\epsilon}-v_{1\epsilon})_{t}\xi_{\delta
t}\nonumber\\
&&\quad-(\alpha
(u_{\epsilon}-v_{1\epsilon})+H_{\epsilon}(u_{\epsilon})-H_{\epsilon}(v_{1\epsilon}))\xi_{\delta
t}\Big]dxdydt \nonumber\\
&&\quad+\int_{\Omega}(\alpha
(u_{\epsilon}(x,T)-v_{1\epsilon})+H_{\epsilon}(u_{\epsilon}(x,T))-H_{\epsilon}(v_{1\epsilon}))\xi_{\delta
}(x,T)dx\nonumber\\
&&\quad=\int_{\Omega}(\alpha(u_{0\epsilon}-v_{1\epsilon})+\chi_{0\epsilon}-H_{\epsilon}(v_{1\epsilon}))\xi_{\delta}(x,0)dx\leq0.
\end{eqnarray}
By Lemma 2.1 of \cite{[Z]}, we obtain from (2.14)
\begin{eqnarray}\label{e2.11 }
&& u_{\epsilon}\leq v_{1\epsilon}\;\; \text{ a.e. in } Q
\end{eqnarray}
and by the monotonicity of $H_{\epsilon}$, we get
\begin{eqnarray}\label{e2.12}
&& H_{\epsilon}(u_{\epsilon})\leq H_{\epsilon}(v_{1\epsilon})\;\; \text{ a.e. in } Q.
\end{eqnarray}

\n We recall that from the proof of existence (see \cite{[G]} or \cite{[Z]} for example), we know that we have up to
a subsequence
\begin{eqnarray}\label{e2.17-2.18}
u_{\epsilon}\rightharpoonup u  \;\;\; &&\text{weakly in} \; L^{2}(0,T; H^{1}(\Omega)).\\
H_{\epsilon}(u_{\epsilon})\rightharpoonup \chi  \;\;\;
&&\text{weakly in} \; L^{2}(Q)
\end{eqnarray}
where $(u,\chi)$ is a solution of problem $(P)$.

\n Similarly, we have since the solution of problem $(P^s_1)$ is unique
\begin{eqnarray}\label{e2.19-2.20}
v_{1\epsilon}\rightharpoonup v_1  \;\;\; &&\text{weakly in} \; H^{1}(\Omega).\\
H_{\epsilon}(v_{1\epsilon})\rightharpoonup \gamma_1  \;\;\;
&&\text{weakly in} \; L^{2}(\Omega).
\end{eqnarray}

\n Now, let $\xi\in\mathcal{D}(Q)$ with $\xi\geq0.$ Passing to the limit, we obtain by using
(2.15)-(2.20) 
\begin{eqnarray*}
&&\int_{Q}(v_{1}-u)\xi dxdt=\lim_{\epsilon\rightarrow0}\int_{Q}(v_{1\epsilon}-u_{\epsilon})\xi
dxdt\geq0,\nonumber\\
&&\int_{Q}(\gamma_1-\chi)\xi dxdt
=\lim_{\epsilon\rightarrow0}\int_{Q}(H_{\epsilon}(v_{1\epsilon})-H_{\epsilon}(u_{\epsilon}))\xi
dxdt\geq0,\nonumber
\end{eqnarray*}
which leads to
\begin{eqnarray}\label{e2.21-2.22}
&& u\leq v_{1}\;\; \text{ a.e. in } Q\\
&& \chi\leq \gamma_1 \;\; \text{ a.e. in } Q.
\end{eqnarray}

\n Similarly, for $\delta>0$ the function
$\displaystyle{\xi_{\delta}=\frac{(v_{0\epsilon}-u_{\epsilon}-\delta)^{+}}{v_{0\epsilon}-u_{\epsilon}}}$
belongs to $ H^1(Q)$ and satisfies $ \xi_{\delta}=0$ on $\Sigma_{2}$ since $\phi_0\leq\phi$ on
$\Sigma_2.$ Then by taking into account (2.9)-(2.10), we get 
\begin{eqnarray}\label{e2.23}
&&\int_{Q}\Big[a(x)\big(\nabla (v_{0\epsilon}-u_{\epsilon})
+(H_{\epsilon}(v_{0\epsilon})-H_{\epsilon}(u_{\epsilon}))e\big).\nabla\xi_{\delta}+\epsilon(v_{0\epsilon}-u_{\epsilon})_{t}\xi_{\delta
t}\nonumber\\
&&\quad-(\alpha
(v_{0\epsilon}-u_{\epsilon})+H_{\epsilon}(v_{0\epsilon})-H_{\epsilon}(u_{\epsilon}))\xi_{\delta
t}\Big]dxdt \nonumber\\
&&\quad+\int_{\Omega}(\alpha
(v_{0\epsilon}-u_{\epsilon}(x,T))+H_{\epsilon}(v_{0\epsilon})-H_{\epsilon}(u_{\epsilon}(x,T)))\xi_{\delta
}(x,T)dx\nonumber\\
&&\quad=\int_{\Omega}(\alpha(v_{0\epsilon}-u_{0\epsilon})+H_{\epsilon}(v_{0\epsilon})-\chi_{0\epsilon})\xi_{\delta}(x,0)dx\leq 0.
\end{eqnarray}

\n By Lemma 2.1 of \cite{[Z]}, we obtain from (2.23)
\begin{eqnarray}\label{e2.24}
&& v_{0\epsilon}\leq u_{\epsilon}\;\; \text{ a.e. in } Q
\end{eqnarray}
and by the monotonicity of $H_{\epsilon}$, we get
\begin{eqnarray}\label{e2.25}
&& H_{\epsilon}(v_{0\epsilon})\leq H_{\epsilon}(u_{\epsilon})\;\; \text{ a.e. in } Q.
\end{eqnarray}

\n Arguing as above and using (2.24)-(2.25), we obtain by passing to the limit
up to a subsequence, that we have for any $\xi\in\mathcal{D}(Q)$ with $\xi\geq0$
\begin{eqnarray*}
&&\int_{Q}(u-v_0)\xi dxdt=\lim_{\epsilon\rightarrow0}\int_{Q}(u_{\epsilon}-v_{0\epsilon})\xi
dxdt\geq0,\nonumber\\
&&\int_{Q}(\chi-\gamma_0)\xi dxdt
=\lim_{\epsilon\rightarrow0}\int_{Q}(H_{\epsilon}(u_{\epsilon})-H_{\epsilon}(v_{0\epsilon}))\xi
dxdt\geq0,\nonumber
\end{eqnarray*}
which leads to
\begin{eqnarray}\label{e2.26-2.27}
&& v_0\leq u\;\; \text{ a.e. in } Q\\
&& \gamma_0\leq \chi \;\; \text{ a.e. in } Q.
\end{eqnarray}

\n Combining (2.21)-(2.22) and (2.26)-(2.27), we obtain (2.11)-(2.12).
\qed

\begin{remark}\label{r2.1} Assume that $a(x)e\in C^{0,1}(\overline{\Omega})$.
Then we get from (2.11)-(2.12) taking into account (2.5) 
\begin{eqnarray}\label{e2.28-2.29}
&&u(x,t)>0  \;\; \text{ if } \; 0<x_2<\epsilon_{0}\\
&&u(x,t)=\chi(x,t)=0 \;\; \text{ if } \; K-\epsilon_{0}<x_2<K.
\end{eqnarray}
\end{remark}

\section{Uniqueness of the solution in rectangular dams}\label{s3}

\vs 0.3cm\n In this section we assume that
\begin{eqnarray}\label{e3.1-2}
&& a(x)\in C^{0,1}(\overline{\Omega}),\quad \text{with}\quad N=\sup_{i,j,k}|(a_{ij})_{x_k}|_{\infty}.\\
&& a(x) \text{ is a symmetric matrix}.
\end{eqnarray}

\vs 0.3cm\n Here is our main result.

\begin{theorem}\label{t3.1} Assume that (2.2)  and (3.1)-(3.2) hold.
Then the solution of the problem $(P)$ associated with the initial
data $(u_0,\chi_0)$ and satisfying (2.28)-(2.29) is unique.
\end{theorem}

\n Let $(u_{1},\chi_{1})$ and $(u_{2},\chi_{2})$ be
two solutions of the problem (P) satisfying (2.28)-(2.29). Set
\begin{eqnarray}
&&w=u_{1}-u_{2}   \;\; \text{ and } \;\; \eta=\alpha
w+\chi_{1}-\chi_{2}.\nonumber
\end{eqnarray}

\vs 0.2cm\n We consider the following problem
\begin{eqnarray}\label{e3.3-3.5}
&&~\text{Find}~v\in L^2(0,T;H^1(\Omega))~\text{ such that}:\nonumber\\
&&~div(a(x)\nabla v)=-\eta\; \text{ in } \; \Omega \; \text{ for each} \; t\in [0,T]\\
&&~v=0  \; \text{ on } \; \Gamma_{2} \\
&&~a(x)\nabla v.\nu=0  \; \text{ on } \; \Gamma_{1}.
\end{eqnarray}

\n The we have
\begin{lemma}\label{l3.1} There exists a unique weak solution of the problem (3.3)-(3.5).
\end{lemma}

\n \emph{Proof.} First, we observe (see \cite{[LZ]}) that 
$\alpha u_i+\chi_i\in C^0([0,T];L^2(\Omega))$, $i=0, 1$.
As a consequence, we have $\eta\in C^0([0,T];L^2(\Omega))$.
Let $V=\{v\in H^1(\Omega)~/~v=0\text{ on }\Gamma_2~\}$.
Then $V$ is a Hilbert space under the $H^1(\Omega)$ norm, and
by applying Lax-Milgram's Theorem and taking into account (1.2)-(1.3), there exists for 
each $t\in[0,T]$ a unique solution $v(x,t)$ 
of the following problem
\begin{eqnarray}\label{e3.6} \left\{
\begin{aligned}
&v(.,t)\in V\\
&\int_\Omega a(x)\nabla v(x,t).\nabla \xi dx=\int_\Omega \eta(x,t) \xi dx \;\;\;\; \forall \xi\in V.
\end{aligned} \right.
\end{eqnarray}
Choosing $\xi\in\mathcal{D}(\Omega)$ in (3.6), we obtain (3.3) in $\mathcal{D}'(\Omega)$
and therefore in $C^0([0,T];L^2(\Omega))$.
(3.4) is satisfied in the trace sense and (3.5). Writing (3.6) for $\xi\in C^\infty(\overline{\Omega})$ with
$\xi=0 $   on $\Gamma_{2}$, and taking into account (3.3), we obtain (3.4) in $H^{-1/2}(\Gamma_1)$. 

\n Choosing $v$ as a test function in (3.6) and using (1.2), H\"{o}lder and Poincar\'{e}'s inequalities,
we obtain
\begin{eqnarray}\label{e3.7}
\int_\Omega |\nabla v(x,t)|^2 dx\leq{1\over\lambda^2}\int_\Omega |\eta(x,t)|^2dx.
\end{eqnarray}
Integrating (3.7) over the interval $[0,T]$ and using the fact that $\chi_i\in L^\infty(Q)$ , $u_i\in L^\infty(0,T;L^\infty(\Omega))$ (see \cite{[LZ]}), we obtain
\begin{eqnarray*}
&&\int_Q |\nabla v(x,t)|^2 dxdt\leq{1\over\lambda^2}\int_0^T\int_\Omega |\eta(x,t)|^2dxdt\leq {{|\eta|_{L^\infty(0,T;L^\infty(\Omega))}}\over\lambda^2}.
\end{eqnarray*}
Using Poincar\'{e}'s inequality, we obtain $v\in L^2(0,T;H^1(\Omega))$. Hence $v$ is the unique solution of (3.3)-(3.5).\\
\qed

\begin{remark}\label{r3.1} 
By the regularity theory (see \cite{[GT]} for example), the solution $v$ of the problem (3.3)-(3.5) satisfies 
$v\in L^2(0,T;C^1(\Omega\cup\mathring{\Gamma}_1\cup\mathring{\Gamma}_2))\cap L^2(0,T;H^2(\Omega\cup\mathring{\Gamma}_1\cup\mathring{\Gamma}_2)).$
\end{remark}

\vs0.3cm\n Now, let us denote by $\widetilde{g}$ the mean with respect to $t$ of a
function $g(x,t)$ defined by
\begin{eqnarray}
\widetilde{g}(x,t)=\frac{1}{h}\int_{t}^{t+h}g(x,s)ds.\nonumber
\end{eqnarray}
Then we have
\begin{eqnarray}\label{e3.8}
&& \widetilde{g}\rightarrow g \;\; \text{ as } \;
h\rightarrow0\nonumber\\
&& \frac{\partial \widetilde{g}}{\partial
t}=\frac{1}{h}(g(x,t+h)-g(x,t)).
\end{eqnarray}
Moreover, it is easy to check that
\begin{eqnarray}\label{e3.9-3.11}
&&~div(a(x)\nabla \widetilde{v})=-\widetilde{\eta}\; \text{ in } \; \Omega \; \text{ for  all } \; t\in [0,T]\\
&&~\widetilde{v}=0  \; \text{ on } \; \Gamma_2 \\
&&~a(x)\nabla \widetilde{v}.\nu=0  \; \text{ on } \; \Gamma_1.
\end{eqnarray}

\n Since $\chi_i=1$ $(i=1,2)$ in a neighborhood of $\Gamma_1$, we obtain from $(P)iii)$
\begin{eqnarray*}
&& (\alpha u_i+\chi_i)_{t}=div(a(x)(\nabla u_i+ \chi_i e ))
~~\text{ in } \mathcal{D}'(Q),\\
&& a(x)\big(\nabla u_i +e\big).\nu=0  \; \text{ on } \; \Gamma_1,~~i=1,2.
\end{eqnarray*}
Using the fact that $u_1=u_2$ on $\Sigma_2, $ and writing the previous two
equations for $(u_1,\chi_1)$ and $(u_2,\chi_2)$
and subtracting them from each other, we get
\begin{eqnarray}\label{e3.12-14}
&& (\alpha w+\chi_1-\chi_2)_{t}=div\big(a(x)(\nabla
w+(\chi_1-\chi_2)e)\big) \;\;\;\;\text{ in }
\mathcal{D}'(Q).\\
&& w=0  \; \text{ on } \; \Sigma_{2}\\
&&a(x)(\nabla w).\nu=0  \; \text{ on } \; \Sigma_1.
\end{eqnarray}

\n Then we have

\begin{lemma}\label{l3.2}
For $h$ small enough we have
\begin{eqnarray}\label{e3.15-16}
&&~\widetilde{\eta}_{t}=div(a(x)(\nabla \widetilde{w}+ (\widetilde{\chi_1}-\widetilde{\chi_2})e) ) \;\;\;\;\text{ in } \mathcal{D}'(Q)\\
&&~\widetilde{w}=0  \; \text{ on } \; \Sigma_2\\
&&~a(x)\nabla \widetilde{w}.\nu=0  \; \text{ on } \; \Sigma_1.
\end{eqnarray}
\end{lemma}
\n \emph{Proof.} (3.16) and (3.17) are a direct consequence of (3.13) and (3.14).
To establish (3.15), let $\zeta\in\mathcal{D}(Q)$ such that for some $\tau_0>0$, $supp(\zeta)\subset
\Omega\times (\tau_0,T-\tau_0)$. We denote by $\widehat{\zeta}$ the
function defined by 
$\displaystyle{\widehat{\zeta}(x,t)=\frac{1}{h}\int_{t-h}^{t}\zeta(x,s)ds}$.
Since for $|h|<\tau_0/2$, the functions $\pm\widehat{\zeta}$ are test functions for problem $(P)$,
we obtain for $i=1, 2$
\begin{eqnarray}\label{e3.18}
&&\int_Q a(x)\big(\nabla u_i +\chi_i e\big)
  \cdot\nabla\widehat{\zeta}dxdt=\int_Q (\alpha u_i+\chi_i)\widehat{\zeta}_t dxdt.
\end{eqnarray}
\n For the right hand side of (3.18), we have by using change of variables
\begin{eqnarray}\label{e3.19}
&&\int_Q (\alpha u_i+\chi_i)\widehat{\zeta}_t dxdt
=\int_Q (\alpha u_i+\chi_i)\frac{1}{h}\zeta(x,t) dxdt-\int_Q (\alpha u_i+\chi_i)\frac{1}{h}\zeta(x,t-h) dxdt\nonumber\\
&&~=\int_Q (\alpha u_i+\chi_i)(x,t)\frac{1}{h}\zeta(x,t) dxdt
-\int_\Omega\int_{-h}^{T-h} \frac{1}{h}(\alpha u_i(x,t+h)+\chi_i(x,t+h))\zeta(x,t)dxdt\nonumber\\
&&~=\int_Q (\alpha u_i+\chi_i)(x,t)\frac{1}{h}\zeta(x,t) dxdt
-\int_\Omega\int_0^T \frac{1}{h}(\alpha u_i(x,t+h)+\chi_i(x,t+h))\zeta(x,t)dxdt\nonumber\\
&&~=-\int_Q \frac{1}{h}[(\alpha u_i(x,t+h)+\chi_i(x,t+h))-(\alpha u_i+\chi_i)(x,t)]\zeta(x,t) dxdt\nonumber\\
&&~=-\int_Q (\widetilde{\alpha u_i+\chi_i})_t\zeta(x,t) dxdt.
\end{eqnarray}

\n For the left hand side of (3.18), we have by integrating by parts 
\begin{eqnarray}\label{e3.20}
&&\int_Q a(x)\big(\nabla u_i +\chi_i e\big)
  \cdot\nabla\widehat{\zeta}dxdt=-\int_Q \Big[\int_0^t a(x)\big(\nabla u_i +\chi_i e\big)ds\Big]\cdot\nabla\widehat{\zeta}_t dxdt\nonumber\\
&&\qquad=-\int_Q \Big[\int_0^t a(x)\big(\nabla u_i +\chi_i e\big)ds\Big]\cdot\nabla\big(\frac{1}{h}(\zeta(x,t)-\zeta(x,t-h))\big) dxdt\nonumber\\
&&\qquad=-\int_Q \frac{1}{h}\Big[\int_0^t a(x)\big(\nabla u_i +\chi_i e\big)ds\Big]\cdot\nabla\zeta(x,t) dxdt\nonumber\\
&&\qquad\quad+\int_Q \frac{1}{h}\Big[\int_0^t a(x)\big(\nabla u_i +\chi_i e\big)ds\Big]\cdot\nabla\zeta(x,t-h) dxdt\nonumber\\
&&\qquad=-\int_Q \frac{1}{h}\Big[\int_0^t a(x)\big(\nabla u_i +\chi_i e\big)ds\Big]\cdot\nabla\zeta(x,t) dxdt\nonumber\\
&&\qquad\quad+\int_\Omega\int_0^T \frac{1}{h}\Big[\int_0^{t+h} a(x)\big(\nabla u_i +\chi_i e\big)ds\Big]\cdot\nabla\zeta(x,t)dxdt\nonumber\\
&&\qquad=\int_Q \frac{1}{h}\Big[\int_t^{t+h} a(x)\big(\nabla u_i +\chi_i e\big)ds\Big]\cdot\nabla\zeta(x,t)  dxdt\nonumber\\
&&\qquad=\int_Q a(x)\big(\nabla \widetilde{u_i} +\widetilde{\chi_i} e\big)\cdot\nabla\zeta(x,t) dxdt.
\end{eqnarray}
\n Then we deduce from (3.18)-(3.20) that
\begin{eqnarray*}
&&\int_Q a(x)\big(\nabla \widetilde{u_i} +\widetilde{\chi_i} e\big)\cdot\nabla\zeta(x,t) dxdt=-\int_Q (\widetilde{\alpha u_i+\chi_i})_t\zeta(x,t) dxdt.
\end{eqnarray*}

\n Writing the last equation for $i=1, 2$ and subtracting the two equations, we get
\begin{eqnarray*}
&&\int_Q a(x)\big(\nabla \widetilde{w} +(\widetilde{\chi_1}-\widetilde{\chi_2}) e\big)\cdot\nabla\zeta(x,t) dxdt=-\int_Q \widetilde{\eta}_t\zeta(x,t) dxdt
\end{eqnarray*}
which is (3.15).
\qed

\vs0.2cm\n To prove Theorem 3.1, we need two more lemmas.

\begin{lemma}\label{l3.3}
\begin{eqnarray}\label{e3.21}
&&\frac{1}{2}\frac{\partial}{\partial
t}\int_{\Omega}a(x)\nabla\widetilde{v}.\nabla\widetilde{v}dx+\int_{\Omega}\widetilde{w}\widetilde{\eta}dx= \int_{\Omega}\alpha
 \widetilde{w}(a_{11}\widetilde{v}_{x_1}+a_{22}\widetilde{v}_{x_2})dx\nonumber\\
&&+\int_{\Omega}(a_{11}\widetilde{v}_{x_1}+a_{12}\widetilde{v}_{x_2})_{x_1}(a_{12}\widetilde{v}_{x_1}+a_{22}\widetilde{v}_{x_2})dx
+\int_{\Omega}(a_{12}\widetilde{v}_{x_1}+a_{22}\widetilde{v}_{x_2})_{x_2}(a_{12}\widetilde{v}_{x_1}+a_{22}\widetilde{v}_{x_2})dx.\nonumber\\
\end{eqnarray}
\end{lemma}
\n \emph{Proof.} From (3.9) and (3.15) we derive
\begin{eqnarray}\label{e3.22}
&&- div(a(x)\nabla \widetilde{v})_{t}=div(a(x)(\nabla
\widetilde{w}+ (\widetilde{\chi_1}-\widetilde{\chi_2})e) )
\;\;\text{ in } \mathcal{D}'(\Omega).
\end{eqnarray}

\n Using (3.10)-(3.11), and taking into account (3.2), we obtain 
\begin{eqnarray}\label{e3.23}
&& <- div(a(x)(\nabla \widetilde{v})_{t},\widetilde{v}>=\int_{\Omega}a(x)
\nabla\widetilde{v}_{t}. \nabla\widetilde{v}dx-\int_{\partial\Omega}a(x)\nabla\widetilde{v}_{t}.\nu\;\widetilde{v}d\sigma(x) \nonumber\\
&&\quad=\int_{\Omega}a(x)\nabla\widetilde{v}_{t}. \nabla\widetilde{v}dx
=\frac{1}{2}\frac{\partial}{\partial t}\int_{\Omega}a(x)\nabla\widetilde{v}.\nabla\widetilde{v}dx.
\end{eqnarray}
\n Similarly, we get by using (3.10), (3.17), and taking into account (3.2)
\begin{eqnarray}
 < div(a(x)\nabla
\widetilde{w}),\widetilde{v}>&=&-\int_{\Omega}a(x)
\nabla\widetilde{w}. \nabla\widetilde{v}dx +\int_{\partial\Omega}a(x)\nabla\widetilde{w}.\nu
\;\widetilde{v}d\sigma(x)\nonumber\\
&=&-\int_{\Omega}a(x) \nabla\widetilde{w}.
\nabla\widetilde{v}dx=-\int_{\Omega}a(x) \nabla\widetilde{v}.
\nabla\widetilde{w}dx \nonumber
\end{eqnarray}
which can be written using (3.9) as
\begin{eqnarray}\label{e3.24}
 < div(a(x)\nabla
\widetilde{w}),\widetilde{v}>&=&-\int_{\Omega}\widetilde{w}\widetilde{\eta}dx.
\end{eqnarray}

\n Similarly, using (3.10) and the fact that $u_1$, $u_2$ satisfy (3.28)-(3.29), 
we get for $\Omega_0=(0,L)\times(\epsilon_0,K-\epsilon_0)$
\begin{eqnarray}\label{e3.25}
&& < div((\widetilde{\chi_{1}}-\widetilde{\chi_{2}})a(x)e),\widetilde{v}>
=-\int_{\Omega_0}(\widetilde{\chi_1}-\widetilde{\chi_2})a(x)e.\nabla\widetilde{v}dx\nonumber\\
&&\quad=\int_{\Omega_0}(-\widetilde{\eta}+\alpha
\widetilde{w})a(x)e.\nabla\widetilde{v}dx\nonumber\\
&&\quad=\int_{\Omega_0}(div(a(x)\nabla\widetilde{v})+\alpha
\widetilde{w})a(x)e.\nabla\widetilde{v}dx\nonumber\\
&&\quad=\int_{\Omega_0}div(a(x)\nabla\widetilde{v}) )a(x)e.\nabla\widetilde{v} dx
+\int_{\Omega_0}\alpha\widetilde{w}a(x)e.\nabla\widetilde{v}dx\nonumber\\
&&\quad=\int_{\Omega_0}(a_{11}\widetilde{v}_{x_1}+a_{12}\widetilde{v}_{x_2})_{x_1}(a_{12}\widetilde{v}_{x_1}+a_{22}\widetilde{v}_{x_2})dx\nonumber\\
&&\quad+\int_{\Omega_0}(a_{12}\widetilde{v}_{x_1}+a_{22}\widetilde{v}_{x_2})_{x_2}(a_{12}\widetilde{v}_{x_1}+a_{22}\widetilde{v}_{x_2})dx\nonumber\\
&&\quad+\int_{\Omega_0}\alpha\widetilde{w}(a_{12}\widetilde{v}_{x_1}+a_{22}\widetilde{v}_{x_2})dx.
\end{eqnarray}
Hence by combining (3.22)-(3.25), we get (3.21).
\qed

\begin{lemma}\label{l3.4}
There exists a positive constant $C$ such that
\begin{eqnarray}\label{e3.26}
&&\frac{\partial}{\partial t}\int_{\Omega}a(x)\nabla\widetilde{v}.\nabla\widetilde{v}dx
+2\int_{\Omega}\widetilde{w}(\widetilde{\chi_{1}}-\widetilde{\chi_{2}})dx\leq
C\int_{\Omega}a(x)\nabla\widetilde{v}.\nabla\widetilde{v}dx.\nonumber\\
&&
\end{eqnarray}
\end{lemma}

\n The proof of Lemma 3.4 requires a lemma.

\begin{lemma}\label{l3.5}
There exists a positive constant $C$ such that
\begin{eqnarray}\label{e3.27-3.30}
&&\int_{\epsilon_0}^{K-\epsilon_0}|\nabla\widetilde{v}(0,x_2)|^2dx_2\leq C\int_{\Omega}|\nabla\widetilde{v}|^2dx\\
&&\int_{\epsilon_0}^{K-\epsilon_0}|\nabla\widetilde{v}(L,x_2)|^2dx_2\leq C\int_{\Omega}|\nabla\widetilde{v}|^2dx\\
&&\int_0^L|\nabla\widetilde{v}(x_1,\epsilon_0)|^2dx_1\leq C\int_{\Omega}|\nabla\widetilde{v}|^2dx\\
&&\int_0^L|\nabla\widetilde{v}(x_1,K-\epsilon_0)|^2dx_1\leq C\int_{\Omega}|\nabla\widetilde{v}|^2dx.
\end{eqnarray}
\end{lemma}
\n \emph{Proof.} $i)$ \n Since $\widetilde{v}\in C^1(\overline{\Omega}_0)$, we have
$\displaystyle{\int_{\epsilon_0}^{K-\epsilon_0}|\nabla\widetilde{v}(0,x_2)|^2 dx_2=\lim_{\delta\rightarrow0}{1\over\delta}\int_0^\delta\int_{\epsilon_0}^{K-\epsilon_0}|\nabla\widetilde{v}|^2dx}$.
It follows that we have for $\delta_1>0$ small  enough
\begin{equation*}
\int_{\epsilon_0}^{K-\epsilon_0}|\nabla\widetilde{v}(0,x_2)|^2dx_2
\leq{2\over{\delta_1}}\int_0^{\delta_1}\int_{\epsilon_0}^{K-\epsilon_0}|\nabla\widetilde{v}|^2dx
\end{equation*}
which gives (3.27).

\vs 0.2cm\n In the same way we establish (3.28), (3.29), and (3.30).
\qed

\n \emph{Proof of Lemma 3.4.} We shall estimate the three integrals in the right hand side of
(3.21). First, we obtain by applying Young's inequality and using (1.2)-(1.3)
\begin{eqnarray}\label{e3.31}
\int_{\Omega_0}\alpha\widetilde{w}(a_{12}\widetilde{v}_{x_1}+
a_{22}\widetilde{v}_{x_2})dx &\leq&\frac{\alpha}{2}\int_{\Omega_0}
 \widetilde{w}^2dx+\frac{\alpha}{2}\int_{\Omega_0}
(a_{12}\widetilde{v}_{x_1}+a_{22}\widetilde{v}_{x_2})^2dx\nonumber\\
&\leq&\frac{\alpha}{2}\int_{\Omega}
 \widetilde{w}^2dx+\frac{\alpha}{2}\int_{\Omega}
(a_{12}^2+a_{22}^2)|\nabla\widetilde{v}|^2dx\nonumber\\
&\leq&\frac{\alpha}{2}\int_{\Omega}\widetilde{w}^{ \;2}dx
+{\alpha \Lambda^2\over\lambda}\int_{\Omega}a(x)\nabla\widetilde{v}.\nabla\widetilde{v}dx.
\end{eqnarray}

\n Next, we have
\begin{eqnarray}\label{e3.32}
&&\int_{\Omega_0}(a_{11}\widetilde{v}_{x_1}+a_{12}\widetilde{v}_{x_2})_{x_1}(a_{12}\widetilde{v}_{x_1}+a_{22}\widetilde{v}_{x_2})dx
=\int_{\Omega_0}(a_{11}\widetilde{v}_{x_1})_{x_1}(a_{12}\widetilde{v}_{x_1})dx\nonumber\\
&&\quad+\int_{\Omega_0}(a_{11}\widetilde{v}_{x_1})_{x_1}(a_{22}\widetilde{v}_{x_2})dx
+\int_{\Omega_0}(a_{12}\widetilde{v}_{x_2})_{x_1}(a_{12}\widetilde{v}_{x_1})dx
+\int_{\Omega_0}(a_{12}\widetilde{v}_{x_2})_{x_1}(a_{22}\widetilde{v}_{x_2})dx\nonumber\\
&&\quad=I_1+I_2+I_3+I_4.
\end{eqnarray}

\n Let us estimate the integrals $I_i$. Expanding and integrating by parts, and
using (1.2)-(1.3), (3.1), and (3.28), we obtain for a positive constant $C_1$
\begin{eqnarray}
I_1&=&\int_{\Omega_0}(a_{11})_{x_1}a_{12}\widetilde{v}_{x_1}^2dx+{1\over2}\int_{\Omega_0}a_{11}a_{12}(\widetilde{v}_{x_1}^2)_{x_1}dx\nonumber\\
&=&\int_{\Omega_0}(a_{11})_{x_1}a_{12}\widetilde{v}_{x_1}^2dx-{1\over2}\int_{\Omega_0}(a_{11}a_{12})_{x_1}\widetilde{v}_{x_1}^2dx\nonumber\\
&&+{1\over2}\int_{\epsilon_0}^{K-\epsilon_0}(a_{11}a_{12}\widetilde{v}_{x_1}^2)(L,x_2)dx_2
-{1\over2}\int_{\epsilon_0}^{K-\epsilon_0}(a_{11}a_{12}\widetilde{v}_{x_1}^2)(0,x_2)dx_2\nonumber
\end{eqnarray}
\begin{eqnarray}\label{e3.33}
&\leq&{1\over2}\int_{\Omega_0}(a_{11})_{x_1}a_{12}\widetilde{v}_{x_1}^2dx
-{1\over2}\int_{\Omega_0}a_{11}(a_{12})_{x_1}\widetilde{v}_{x_1}^2dx
+{{\Lambda^2}\over2}\int_{\epsilon_0}^{K-\epsilon_0}\widetilde{v}_{x_1}^2(L,x_2)dx_2\nonumber\\
&\leq&\Lambda N\int_{\Omega_0}\widetilde{v}_{x_1}^2dx
+{{\Lambda^2}\over2}\int_{\epsilon_0}^{K-\epsilon_0}\widetilde{v}_{x_1}^2(L,x_2)dx_2
\leq C_1\int_{\Omega}a(x)\nabla\widetilde{v}.\nabla\widetilde{v}dx.
\end{eqnarray}

\n Note that since $\widetilde{v}(0,x_2)=\widetilde{v}(L,x_2)=0$ for $0<x_2<K$, we have
\begin{equation}\label{3.34}
\widetilde{v}_{x_2}(0,x_2)=\widetilde{v}_{x_2}(L,x_2)=0\quad \text{for}\quad
0<x_2<K.
\end{equation}
Integrating by parts and using (3.34), we obtain
\begin{eqnarray}\label{e3.35}
I_2&=&-\int_{\Omega_0}a_{11}\widetilde{v}_{x_1}(a_{22}\widetilde{v}_{x_2})_{x_1}dx=
-\int_{\Omega_0}a_{11}(a_{22})_{x_1}\widetilde{v}_{x_1}\widetilde{v}_{x_2}dx
-{1\over2}\int_{\Omega_0}a_{11}a_{22}(\widetilde{v}_{x_1}^2)_{x_2}dx\nonumber\\
&=&-\int_{\Omega_0}a_{11}(a_{22})_{x_1}\widetilde{v}_{x_1}\widetilde{v}_{x_2}dx
+{1\over2}\int_{\Omega_0}a_{11}(a_{22})_{x_2}\widetilde{v}_{x_1}^2dx
+{1\over2}\int_{\Omega_0}(a_{11})_{x_2} a_{22}\widetilde{v}_{x_1}^2dx\nonumber\\
&&-{1\over2}\int_0^L (a_{11}a_{22}\widetilde{v}_{x_1}^2)(x_1,K-\epsilon_0)dx_1
+{1\over2}\int_0^L (a_{11}a_{22}\widetilde{v}_{x_1}^2)(x_1,\epsilon_0)dx_1.
\end{eqnarray}
\n Using Young's inequality, (1.2)-(1.3), (3.1), and (3.29), we obtain from (3.35), for a 
positive constant $C_2$ 
\begin{eqnarray}\label{e3.36}
I_2&\leq& \Lambda N\int_{\Omega_0}|\widetilde{v}_{x_1}|.|\widetilde{v}_{x_2}|dx
+\Lambda N\int_{\Omega_0}\widetilde{v}_{x_1}^2dx
+{{\Lambda^2}\over2}\int_0^L \widetilde{v}_{x_1}^2(x_1,\epsilon_0)dx_1\nonumber\\
&\leq& {{3\Lambda N}\over2}\int_{\Omega_0}|\nabla\widetilde{v}|^2dx
+{{\Lambda^2}\over2}\int_0^L \widetilde{v}_x^2(x_1,\epsilon_0)dx_1
\leq C_2\int_{\Omega}a(x)\nabla\widetilde{v}.\nabla\widetilde{v}dx.
\end{eqnarray}

\n Expanding and integrating by parts, and using Young's inequality, (1.2)-(1.3), (3.1),
and (3.30), we obtain for a positive constant $C_3$ 
\begin{eqnarray}\label{e3.37}
&&I_3=\int_{\Omega_0}(a_{12})_{x_1}a_{12}\widetilde{v}_{x_1}\widetilde{v}_{x_2}dx
+{1\over2}\int_{\Omega_0}a_{12}^2(\widetilde{v}_{x_1}^2)_{x_2}dx\nonumber\\
&&~~=\int_{\Omega_0}(a_{12})_{x_1}a_{12}\widetilde{v}_{x_1}\widetilde{v}_{x_2}dx
-{1\over2}\int_{\Omega_0}(a_{12}^2)_{x_2}\widetilde{v}_{x_1}^2dx\nonumber\\
&&~~~+{1\over2}\int_0^L (a_{12}^2\widetilde{v}_{x_1}^2)(x_1,K-\epsilon_0)dx_1
-{1\over2}\int_0^L (a_{12}^2\widetilde{v}_{x_1}^2)(x_1,\epsilon_0)dx_1\nonumber\\
&&~~\leq {\Lambda N\over2}\int_{\Omega_0}(\widetilde{v}_{x_1}^2+\widetilde{v}_{x_2}^2)dx
-{1\over2}\int_{\Omega_0}(a_{12}^2)_{x_2}\widetilde{v}_{x_1}^2dx
+{{\Lambda^2}\over2}\int_0^L \widetilde{v}_{x_1}^2(x_1,K-\epsilon_0)dx_1\nonumber\\
&&~~\leq C_3\int_{\Omega}a(x)\nabla\widetilde{v}.\nabla\widetilde{v}dx.
\end{eqnarray}
\n Integrating by parts and using (3.34), and taking into account (1.2)-(1.3) and (3.1),
we obtain for a positive constant $C_4$
\begin{eqnarray}\label{e3.38}
I_4&=&\int_{\Omega_0}(a_{12})_{x_1}a_{22}\widetilde{v}_{x_2}^2dx
+{1\over2}\int_{\Omega_0}a_{12}a_{22}(\widetilde{v}_{x_2}^2)_{x_1}dx\nonumber\\
&=&\int_{\Omega_0}(a_{12})_{x_1}a_{22}\widetilde{v}_{x_2}^2dx
-{1\over2}\int_{\Omega_0}(a_{12}a_{22})_{x_1}\widetilde{v}_{x_2}^2dx\nonumber\\
&=&{1\over2}\int_{\Omega_0}\big((a_{12})_{x_1}a_{22}-a_{12}(a_{22})_{x_1}\big)\widetilde{v}_{x_2}^2dx
\leq C_4\int_{\Omega}a(x)\nabla\widetilde{v}.\nabla\widetilde{v}dx.
\end{eqnarray}

\n To estimate the last integral in the right hand side of (3.21), we integrate by parts and 
use (3.29-(3.30) and (1.2). We obtain for a positive constant $C_5$
\begin{eqnarray*}
&&\int_{\Omega_0}(a_{12}\widetilde{v}_{x_1}+a_{22}\widetilde{v}_{x_2})_{x_2}(a_{12}\widetilde{v}_{x_1}+a_{22}\widetilde{v}_{x_2})dx
={1\over2}\int_{\Omega_0}((a_{12}\widetilde{v}_{x_1}+a_{22}\widetilde{v}_{x_2})^2)_{x_2}dx\nonumber\\
&&\quad={1\over2}\int_0^L(a_{12}\widetilde{v}_{x_1}+a_{22}\widetilde{v}_{x_2})^2(x_1,K-\epsilon_0)dx_1
-{1\over2}\int_0^L(a_{12}\widetilde{v}_{x_1}+a_{22}\widetilde{v}_{x_2})^2(x_1,\epsilon_0)dx\nonumber\\
\end{eqnarray*}
\begin{eqnarray}\label{e3.39}
&&\quad\leq{1\over2}\int_0^L(a_{12}\widetilde{v}_{x_1}+a_{22}\widetilde{v}_{x_2})^2(x_1,K-\epsilon_0)dx_1\nonumber\\
&&\quad\leq\int_0^L((a_{12}\widetilde{v}_{x_1})^2+(a_{22}\widetilde{v}_{x_2})^2)(x_1,K-\epsilon_0)dx_1\nonumber\\
&&\quad\leq\Lambda^2\int_0^L((\widetilde{v}_{x_1})^2+(\widetilde{v}_{x_2})^2)(x_1,K-\epsilon_0)dx_1
\leq C_5\int_{\Omega}a(x)\nabla\widetilde{v}.\nabla\widetilde{v}dx.
\end{eqnarray}

\n Finally, combining (3.21), (3.31)-(3.33), and (3.35)-(3.39), we get for a positive constant $C$

\n $\displaystyle{\frac{\partial}{\partial t}\int_{\Omega}a(x)\nabla\widetilde{v}.\nabla\widetilde{v}dx
+2\int_{\Omega}\widetilde{w}(\widetilde{\chi_{1}}-\widetilde{\chi_{2}})dx\leq
C\int_{\Omega}a(x)\nabla\widetilde{v}.\nabla\widetilde{v}dx}$,
which is (3.26).
\qed

\vs 0.3cm\n \emph{Proof of Theorem 3.1.} First, integrating (3.26) from $0$ to $t$ and letting $h\rightarrow0$, we get
\begin{eqnarray}\label{e3.40}
\int_{\Omega}a(x)\nabla v.\nabla v dx
+2\int_{0}^{t}\int_{\Omega}w(\chi_{1}-\chi_{2})dxds
\leq C\int_{0}^{t}\int_{\Omega}a(x)\nabla v.\nabla vdxds.
\end{eqnarray}
\n Next, we observe that since $u_i\in H(\chi_i)$ a.e. in $Q$, we have 
\begin{equation}\label{e3.41}
w(\chi_{1}-\chi_{2})\geq0\quad \text{ a.e. in }~~Q.
\end{equation}
\n Setting $\displaystyle{F(t)= \int_{0}^{t}\int_{\Omega}a(x)\nabla v.\nabla vdxds}$,
we deduce from (3.40)-(3.41) that
\begin{eqnarray}\label{e3.42}
&&F'(t)\leq CF(t)  \;\;\; \forall t\in [0,T].
\end{eqnarray}
\n Integrating (3.42), we get since $F(0)=0$,
$~~ 0\leq F(t)\leq F(0)e^{Ct}=0 \;\;\; \forall t\in [0,T]$, or
\begin{equation*}
\int_{0}^{t}\int_{\Omega}a(x)\nabla v.\nabla vdxds=0\quad \forall t\in [0,T].
\end{equation*}
Using (1.2), we obtain $~\nabla v=0$~ a.e. in ~~$Q$.
Taking into account that $v=0$ on $\Gamma_{2}\subset
\partial\Omega$ and the connectedness of $\Omega$, we obtain $v=0$ in
$Q.$ Going back to (3.3), we obtain $\eta=0$ a.e. in $Q$, which reads
\begin{eqnarray}\label{e3.43}
w+\chi_{1}-\chi_{2}=0  \;\; \text{ in } \; Q.
\end{eqnarray}
Multiplying (3.43) by $w,$ we get $~~\alpha w^{2}+w(\chi_1-\chi_2)=0~~$ a.e. in 
$~~Q$. Taking into account (3.41), we obtain 
$w^2=0$ a.e. in $Q$, or $u_1=u_2$  a.e. in $Q$.
Finally, we obtain from (3.43) that $\chi_1=\chi_2$ a.e. in $Q$.
This achieves the proof.
\qed

\vs 0,3cm
\n\emph{Acknowledgments } The second author is grateful to Prof. J. F. Rodrigues
for kindly inviting him to the CMAF where he enjoyed excellent research conditions
during the preparation of his Ph.D Thesis.

\end{document}